\documentclass[12pt,reqno]{amsart}

\usepackage{pgfplots}
\usepackage{times}
\usepackage{amsmath,amsfonts, amstext,amssymb,amsbsy,amsopn,amsthm}
\usepackage{dsfont}
\usepackage{esint}
\usepackage{graphicx}   
\usepackage{hyperref}
\usepackage[all]{xy}
\usepackage{color}
\usepackage{tikz}
\usepackage{lineno,hyperref}

\newcommand{\fr}{\displaystyle\frac}
\newcommand{\jf}{\displaystyle\int}
\newcommand{\lt}{\left}
\newcommand{\rt}{\right}

\newcommand{\mb}{\mbox}

\newcommand{\tm}{\times}

\newcommand{\pl}{\partial}
\newcommand{\al}{\alpha}

\newcommand{\abs}[1]{\lvert#1\rvert}
\newcommand{\Abs}[1]{\left\lvert#1\right\rvert}

\newcommand{\R}{\mathbb{R}}

\newtheorem{theorem}{Theorem}[section]
\newtheorem{lemma}[theorem]{Lemma}

\theoremstyle{definition}
\newtheorem{definition}[theorem]{Definition}
\theoremstyle{remark}
\newtheorem{remark}{Remark}[section]
\newtheorem{counterexample}{Counterexample}

\theoremstyle{remark}

\numberwithin{equation}{section}
\allowdisplaybreaks[4]

\begin{document}

\title{Liouville theorem for fully fractional master equations and its applications}

\author{Wenxiong Chen }
\address{Department of Mathematical Sciences, Yeshiva University, New York, NY,  10033 USA}
\email{wchen@yu.edu}

\author{Lingwei Ma}
\address{School of Mathematical Sciences, Tianjin Normal University,
Tianjin, 300387, P.R. China, and Department of Mathematical Sciences, Yeshiva University, New York, NY, 10033, USA}
\email{mlw1103@outlook.com}

\author{Yahong Guo}
\address{School of Mathematical Sciences, Nankai University, Tianjin, 300071, P.R. China, and Department of Mathematical Sciences, Yeshiva University, New York, NY, 10033, USA}
\email{1120200036@mail.nankai.edu.cn}

\date{\today}

\begin{abstract}
In this paper, we study the fully fractional master equation
\begin{equation}\label{pdeq1}
  (\partial_t-\Delta)^s u(x,t) =f(x,t,u(x,t)),\,\,(x, t)\in \mathbb{R}^n\times \mathbb{R}.
\end{equation}

First we prove a Liouville type theorem for the homogeneous equation
\begin{equation}\label{pdeq0}
  (\partial_t-\Delta)^s u(x,t) = 0,\,\,(x, t)\in \mathbb{R}^n\times \mathbb{R},
\end{equation}
where $0<s<1$.
When $u$ belongs to the slowly increasing function space
$$\mathcal{L}^{2s,s}(\mathbb{R}^n\times\mathbb{R})=\left\{u(x,t) \in L^1_{\rm loc} (\mathbb{R}^n\times\mathbb{R}) \mid \int_{-\infty}^{+\infty} \int_{\mathbb{R}^n} \frac{|u(x,t)|}{1+|x|^{n+2+2s}+|t|^{\frac{n}{2}+1+s}}\operatorname{d}\!x\operatorname{d}\!t<\infty\right\} $$
and satisfies an additional asymptotic assumption
$$\liminf_{|x|\rightarrow\infty}\frac{u(x,t)}{|x|^\gamma}\geq 0 \; ( \mbox{or} \; \leq 0)  \,\,\mbox{for some} \;0\leq\gamma\leq 1, $$
in the case $\frac{1}{2}<s < 1$,
we prove that all solutions of (\ref{pdeq0}) must be constant. This result includes the previous Liouville theorems on harmonic functions \cite{ABR} and on $s$-harmonic functions \cite{CDL} as special cases.

Then we establish the equivalence between nonhomogeneous pseudo-differential equations (\ref{pdeq1}) and the corresponding integral equations. We believe that these integral equations will become very useful tools in further analysing qualitative properties of solutions, such as regularity, monotonicity, and symmetry.

In the process of deriving the Liouville type theorem,  through very delicate calculations, we obtain an optimal estimate on
the decay rate of $(\partial_t-\Delta)_{\rm right}^s \varphi(x,t)$. This sharp estimate will become a key ingredient and an important tool in investigating master equations.

\bigskip


\bigskip

\end{abstract}

\maketitle

\section{Introduction}\label{sec1}

The classical Liouville theorem states that any bounded harmonic function, defined in the entire space $\mathbb{R}^n$ and satisfying
$\Delta u(x)=0$,
must be identically constant.

This result was initially proposed by Liouville and later proven by Cauchy \cite{Cau} using Cauchy's integral formula.

The bounded-ness condition in the Liouville theorem can be relaxed to either one-side bounded-ness or even weaker growth condition at infinity
\begin{equation*}
  \liminf_{|x|\rightarrow\infty}\frac{u(x)}{|x|}\geq 0.
\end{equation*}
For more details, please refer to the monograph \cite{ABR}.

This Liouville theorem has found extensive applications in the analysis of partial differential equations. It has been instrumental in deriving a priori estimates and establishing the qualitative properties of solutions, including their existence, nonexistence, and uniqueness.

Due to the significant role played by classical Liouville theorem, it has been widely studied and extended to various types of elliptic equations  using  diverse methods, including Harnack inequalities,  blow up and compactness arguments, as well as Fourier analysis (cf. \cite{BKN, CDL, Fa, Mos, SZ} and the references therein).

Among them, we would like to particularly mention the work of \cite{BKN} which shows that the classical Liouville theorem have analogues in the $s$-harmonic setting. More precisely, if
$u$ is a nonnegative solution of the following elliptic fractional equation
\begin{equation*}
  (-\Delta)^su(x)=0\,\,\mbox{in}\,\, \mathbb{R}^n,
\end{equation*}
then $u$ must be constant. Subsequently, Chen, D' Ambrosio and Li \cite{CDL} weaken
this one-sided bounded-ness condition to
\begin{equation*}
  \liminf_{|x|\rightarrow\infty}\frac{u(x)}{|x|^\gamma}\geq 0 \,\,\mbox{for}\,\, 0\leq\gamma\leq\min\{1,2s\}.
\end{equation*}

However, for entire solutions of the parabolic equation
\begin{equation}\label{heat}
  \partial_tu(x,t)-\Delta u(x,t)=0\,\, \mbox{in}\,\, \mathbb{R}^n\times\mathbb{R},
\end{equation}
the one-side bounded-ness condition is no longer sufficient to guarantee the conclusion of the Liouville theorem,
as is shown by the counter example
$$u(x,t)=e^{x_1+x_2+...+x_n+nt},$$
which is a positive non-constant solution of \eqref{heat} in $\mathbb{R}^n\times \mathbb{R}$.

When restricted the domain of $t$ to $(-\infty, 0]$, that is when considering the ancient solutions of (\ref{heat}),
Widder \cite{Wid} proved that all bounded solutions $u(x,t)$ must be constant.
  Hirschman \cite{Hir} relaxed the bounded-ness condition to one-side bounded-ness with the additional growth condition at $t=0$ as follows
  $$\liminf_{r\rightarrow\infty}\frac{\log M(r)}{r}\leq 0,$$
  where
  $$M(r)=\sup\{u(x,0)\mid |x|\leq r\}.$$


For the nonlocal parabolic equation
 \begin{equation*}
     \partial_tu(x,t)+(-\Delta)^s u(x,t)=0, \;\; (x,t) \in \mathbb{R}^n\times(-\infty,0],
  \end{equation*}
Serra \cite{Ser} showed that $u$ is a constant  under the growth condition
$$\|u\|_{L^\infty(Q_R^{2s})}\leq C R^\beta\,\, \mbox{for all}  \,\, R\geq1 \,\,\mbox{and some}\,\, 0<\beta<\min\{1, 2s\},$$
where the parabolic cylinder
$$Q_R^{2s}=\{(x,t)\mid |x|< R \,\,\mbox{and}\,\, -R^{2s}<t<0\}.$$

Recently, Ma, Guo, and Zhang \cite{MGZ} proved that the bounded entire solutions of the homogeneous master equation
\begin{equation}\label{model}
    (\partial_t-\Delta)^s u(x,t)=0 \,\,  \mbox{in}\,\,  \mathbb{R}^n\times\mathbb{R},
\end{equation}
must be constant.
Nevertheless, it remains uncertain whether the Liouville theorem for master equation \eqref{model} still holds when weakening the two-sided bounded-ness condition. This is the primary objective of our research here.

The fully fractional heat operator $(\partial_t-\Delta)^s$ was first proposed by Riesz \cite{Riesz}, which can be defined pointwise by the following singular integral
\begin{equation}\label{nonlocaloper}
(\partial_t-\Delta)^s u(x,t)
:=C_{n,s}\int_{-\infty}^{t}\int_{\mathbb{R}^n}
  \frac{u(x,t)-u(y,\tau)}{(t-\tau)^{\frac{n}{2}+1+s}}e^{-\frac{|x-y|^2}{4(t-\tau)}}\operatorname{d}\!y\operatorname{d}\!\tau,
\end{equation}
where $0<s<1$, the integral in $y$ is taken in the Cauchy principal value sense, and the normalization positive constant $$C_{n,s}=\frac{1}{(4\pi)^{\frac{n}{2}}|\Gamma(-s)|},$$
with $\Gamma(\cdot)$ denoting the Gamma function. The singular integral in \eqref{nonlocaloper} converges for each $(x,t) \in \mathbb{R}^n\times\mathbb{R}$ provided
 $$u(x,t)\in C^{2s+\epsilon,s+\epsilon}_{x,\, t,\, {\rm loc}}(\mathbb{R}^n\times\mathbb{R}) \cap \mathcal{L}(\mathbb{R}^n\times\mathbb{R})$$
for some $\varepsilon>0$, where
the space $\mathcal{L}(\mathbb{R}^n\times\mathbb{R})$ is defined by
$$ \mathcal{L}(\mathbb{R}^n\times\mathbb{R}):=\left\{u(x,t) \in L^1_{\rm loc} (\mathbb{R}^n\times\mathbb{R}) \mid \int_{-\infty}^t \int_{\mathbb{R}^n} \frac{|u(x,\tau)|e^{-\frac{|x|^2}{4(t-\tau)}}}{1+(t-\tau)^{\frac{n}{2}+1+s}}\operatorname{d}\!x\operatorname{d}\!\tau<\infty,\,\, \forall \,t\in\mathbb{R}\right\},$$
and the definition of the local parabolic H\"{o}lder space $C^{2s+\epsilon,s+\epsilon}_{x,\, t,\, {\rm loc}}(\mathbb{R}^n\times\mathbb{R})$ will be specified in Section \ref{2}\,. In particular, if $u$ is bounded, we can ensure the integrability of \eqref{nonlocaloper} by assuming only that $u$ is local parabolic H\"{o}lder continuous.
We notice that the operator $(\partial_t-\Delta)^s$ is nonlocal both in space and time, since the value of $(\partial_t-\Delta)^s u$ at a given point $(x,t)$ depends on the values of $u$ over the whole $\mathbb{R}^n$ and even on all the past time before $t$.
The intriguing aspect of this problem is that applying the space-time nonlocal operator $(\partial_t-\Delta)^s$ to a function that only depends on either space or time, it reduces to a familiar fractional order operator, as discussed in \cite{ST}. More precisely, if $u$ is only a function of $x$, then
 \begin{equation*}
   (\partial_t-\Delta)^s u(x)=(-\Delta)^s u(x),
 \end{equation*}
where $(-\Delta)^s$ is the well-known fractional Laplacian of order $2s$.
In recent decades, the well-posedness of solutions to elliptic equations involving the fractional Laplace operator has been extensively investigated, interested readers can refer to \cite{CLL, CLZ, CLZ1, CWu1, DL, DLL, DLQ, LL, LLW, LW, LZ, MZ2, ZL} and references therein.
While if $u=u(t)$, then
 \begin{equation*}
   (\partial_t-\Delta)^s u(t)=\partial_t^s u(t),
 \end{equation*}
where $\partial_t^s$ is usually denoted by $D_{\rm left}^s$, representing the Marchaud left fractional derivative of order $s$, defined as
\begin{equation*}
D_{\rm left}^s u(t)=\frac{1}{|\Gamma(-s)|}\int_{-\infty}^t \frac{u(t)-u(\tau)}{(t-\tau)^{1+s}}\operatorname{d}\!\tau.
\end{equation*}
Moreover, it should be noted that as $s$ tends to $1$ from the left side,
the fractional power of heat operator $(\partial_t-\Delta)^s$
converges to the local heat operator $\partial_t-\Delta$ (cf. \cite{FNW}).

The master equation has a wide range of applications in physical and biological phenomena, such as anomalous diffusion \cite{KBS}, chaotic dynamics \cite{Z}, biological invasions \cite{BRR}, among others. In addition to these areas, it has also been employed in the financial field \cite{RSM}, where it can model the correlation between waiting times and price jumps in transactions. From a probabilistic perspective, the master equation plays a crucial role in the theory of continuous time random walk, where $u$ represents the distribution of particles that make random jumps simultaneously with random time lags (cf. \cite{MK}).
This is in contrast to the nonlocal parabolic equations
\begin{equation}\label{frac-para}
  \partial_t u+(-\Delta)^s u=f
\end{equation}
or the dual fractional parabolic equation
\begin{equation}\label{frac-para1}
  \partial_t^\alpha u+(-\Delta)^s u=f,
\end{equation}
where jumps are independent of the waiting times.
In other words, the master equation takes into account the strong correlation between the waiting times and the particle jumps, whereas the nonlocal parabolic equation \eqref{frac-para} or the dual fractional parabolic equation \eqref{frac-para1} do not.
It is evident that the master equation is of great importance in various fields, and continuous research on it can drive us towards a deeper understanding of complex phenomena.
In the last decade, significant advancements have been made in the research on the qualitative properties of solutions to this class of master equations, such as regularity, monotonicity and symmetry. For details, please refer to  \cite{ACM, CM1, CS2, FD, MGZ, ST} and the references therein. While in the study of  the other two types of fractional parabolic equations (\ref{frac-para}) and (\ref{frac-para1}), there have seen even more progresses, for instance, see \cite{CM2, CWNH, CWu2, ChenWu, CWW, WuC} and the references therein.

It is worth mentioning that \cite{FD, ST} indicate that the fully fractional heat operator $(\partial_t-\Delta)^s$ has a fundamental solution
\begin{equation*}
  G(x,t)=\frac{1}{(4\pi)^{\frac{n}{2}}\Gamma(s)}\frac{e^{-\frac{|x|^2}{4t}}}
  {t^{\frac{n}{2}+1-s}}\chi_{\{t>0\}}
\end{equation*}
satisfying
\begin{equation*}
 (\partial_t-\Delta)^sG(x,t)=\delta(x,t),
\end{equation*}
where $\chi$ is the characteristic function and $\delta$ is the Dirac Delta function. Furthermore, they proved that
$$u(x,t)= \int_{-\infty}^{+\infty} \int_{\mathbb{R}^n} G(x-y,t-\tau)h(y,\tau)\operatorname{d}\!y\operatorname{d}\!\tau$$
is a solution to the master equation
 \begin{equation*}
   (\partial_t-\Delta)^s u(x,t)=h(x,t),
 \end{equation*}
where $h$ is regular enough to ensure that the above integral representation formula
is well defined.

Then it is natural to ask:  {\em Is this solution unique?}

To answer this question, we establish a Liouville type theorem in this context.
Under appropriate assumptions, we prove that all entire solutions of the homogeneous master equation
\begin{equation}\label{model-Lio}
(\partial_t-\Delta)^s u(x,t) =0, \; \,\,(x, t)\in \mathbb{R}^n\times \mathbb{R}
\end{equation}
must be constant. More precisely, we have

\begin{theorem}\label{Liouville}
Let $0<s<1$ and $n\geq2$.
Assume that
$$u\in\mathcal{L}^{2s,s}(\mathbb{R}^n\times\mathbb{R})=\left\{u(x,t) \in L^1_{\rm loc} (\mathbb{R}^n\times\mathbb{R}) \mid \int_{-\infty}^{+\infty} \int_{\mathbb{R}^n} \frac{|u(x,t)|}{1+|x|^{n+2+2s}+|t|^{\frac{n}{2}+1+s}}\operatorname{d}\!x\operatorname{d}\!t<\infty\right\} $$
is a solution of (\ref{model-Lio}) in the sense of distribution.

In the case $\frac{1}{2}<s < 1$, we assume additionally that
\begin{equation}\label{AA}
  \liminf_{|x|\rightarrow\infty}\frac{u(x,t)}{|x|^\gamma}\geq 0 \; ( \mbox{or} \; \leq 0) \,\,\mbox{for some} \;0\leq\gamma\leq 1.
\end{equation}
Then $u$ must be a constant.
\end{theorem}

We say that $u\in \mathcal{L}^{2s,s}(\mathbb{R}^n\times\mathbb{R})$ is a solution of (\ref{model-Lio}) in the sense of distribution, if
\begin{equation}\label{dissol}
\int_{-\infty}^{+\infty} \int_{\mathbb{R}^n} u(x,t)\overline{(\partial_t-\Delta)_{\rm right}^s\varphi(x,t)}\operatorname{d}\!x\operatorname{d}\!t = 0
\end{equation}
for any $\varphi\in \mathcal{S}(\mathbb{R}^n\times\mathbb{R})$, where
the right fractional derivative
\begin{equation}\label{master-R}
(\partial_t-\Delta)_{\rm right}^s \varphi(x,t)
:=C_{n,s}\int_{t}^{+\infty}\int_{\mathbb{R}^n}
  \frac{\varphi(x,t)-\varphi(y,\tau)}{(\tau-t)^{\frac{n}{2}+1+s}}e^{-\frac{|x-y|^2}{4(\tau-t)}}\operatorname{d}\!y\operatorname{d}\!\tau
\end{equation}
considers all future time points after a given time $t$.
In the subsequent section, we will demonstrate that the space $ \mathcal{L}^{2s,s}$ is appropriate for defining the distributional left fractional derivative $(\partial_t-\Delta)^s$.

\begin{remark}
It is noteworthy that the asymptotic assumption \eqref{AA} of $u$ in Liouville Theorem \ref{Liouville} is necessary. When this condition is violated, there are counter examples, such as the function $u(x,t)=x_1$, which is a nonconstant solution of master equation \eqref{model-Lio} and does not satisfy \eqref{AA}.
\end{remark}

\begin{remark}
In particular, if $u$ is a function of $x$ only, i.e., $u(x,t)=u(x)$, then a straightforward calculation yields that
\[\int_{-\infty}^{+\infty}\int_{\mathbb{R}^n} \frac{|u(x)|}{1+|x|^{n+2+2s}+|t|^{\frac{n}{2}+1+s}}\operatorname{d}\!x\operatorname{d}\!t\sim\int_{\mathbb{R}^n} \frac{|u(x)|}{1+|x|^{n+2s}}\operatorname{d}\!x.\]
This indicates that the space $\mathcal{L}^{2s,s}(\mathbb{R}^n\times\mathbb{R})$ is consistent with our familiar space
$$\mathcal{L}^s=\left\{u(x) \in L^1_{\rm loc} (\mathbb{R}^n) \mid \int_{\mathbb{R}^n} \frac{|u(x)|}{1+|x|^{n+2s}}\operatorname{d}\!x<\infty\right\},$$
which is introduced to ensure the well-definition of the fractional Laplacian $(-\Delta)^s$.
From here, one can see that Theorem \ref{Liouville} here includes the Liouville theorem for $s$-harmonic functions established in \cite{CDL} as a special case.
\end{remark}
\medskip

To prove the Liouville theorem, we employ the Fourier transform on tempered distributions. To this end, it is imperative  to derive a key estimate on
the decay rate of $(\partial_t-\Delta)_{\rm right}^s \varphi(x,t)$. Through very delicate calculations, we obtain
\begin{lemma}\label{lm1}
If $\varphi(x,t)\in \mathcal{S}(\R^{n}\times \R),$ then $(\partial_t-\Delta)_{\rm{right}}^s\varphi(x,t)\in C^{\infty}(\R^n\times \R)$ and
\begin{equation}\label{1}
\Abs{(\partial_t-\Delta)_{{\rm{right}}}^s\varphi(x,t)} \leq\fr{C}{1+\abs{x}^{n+2+2s}+\abs{t}^{\frac{n}{2}+1+s}}\end{equation}
for all $(x,t)\in\R^n\tm\R$.
\end{lemma}
The decay rate provided in (\ref{1}) is optimal, as will be illustrated by a counter example in Section \ref{sec2}.  We believe that this sharp estimate will become a key ingredient and an important tool in the analysis of master equations.
\medskip

As an application of the above Liouville theorem for homogeneous master equation, we establish an equivalence between the fully fractional heat equation and the corresponding integral equation.

\begin{theorem}\label{equivalence}
Let $u(x,t)\in  C^{2s+\epsilon,s+\epsilon}_{x,\, t,\, {\rm loc}}(\mathbb{R}^n\times\mathbb{R}) \cap \mathcal{L}^{2s,s}(\mathbb{R}^n\times\mathbb{R})$ be positive and $f$ be nonnegative and continuous. Assume that
$$f(x,t,\zeta)\geq C_2>0 \;\mbox{ for } \;\zeta\geq C_1>0, \;\mbox{ uniformly in } (x,t)\in \mathbb{R}^n\times \mathbb{R}.$$

Then the pseudo-differential equation
\begin{equation}\label{pdeq}
  (\partial_t-\Delta)^s u(x,t) =f(x,t,u(x,t)),\,\,(x, t)\in \mathbb{R}^n\times \mathbb{R},
\end{equation}
is equivalent to the integral equation
\begin{equation}\label{ieq}
  u(x,t)= \int_{-\infty}^{+\infty} \int_{\mathbb{R}^n} G(x-y,t-\tau)f(y,\tau,u(y,\tau))\operatorname{d}\!y\operatorname{d}\!\tau,
\end{equation}
where
\begin{equation*}
  G(x-y,t-\tau)=A_{n,s}\frac{e^{-\frac{|x-y|^2}{4(t-\tau)}}}
  {(t-\tau)^{\frac{n}{2}+1-s}}\chi_{\{t>\tau\}}
\end{equation*}
with $A_{n,s}:=\frac{1}{(4\pi)^{\frac{n}{2}}\Gamma(s)}$.
\end{theorem}

To conclude this section, we will provide a brief outline of the structure of this paper.
In Section \ref{sec2}, we introduce the parabolic H\"{o}lder spaces, the fractional derivative and the Fourier transform in the sense of distributions, and obtain optimal estimate \eqref{1}. In Section \ref{sec3}, we prove our main results,  Theorem \ref{Liouville} and Theorem \ref{equivalence}.

\section{Preliminaries}\label{sec2}

In this section, we first recall the definition of the parabolic H\"{o}lder space and then establish a useful optimal estimate, thus defining the fractional derivative $(\partial_t-\Delta)^s$ of a distribution. Throughout this paper, we use $C$ to denote a general constant whose value may vary from line to line.

\subsection{Parabolic H\"{o}lder space}\label{2.1}
Now we start by stating the definition of parabolic H\"{o}lder space $C^{2\alpha,\alpha}_{x,\, t}(\mathbb{R}^n\times\mathbb{R})$ (cf. \cite{Kry}) as follows.
\begin{itemize}
\item[(i)]
When $0<\alpha\leq\frac{1}{2}$, if
$u(x,t)\in C^{2\alpha,\alpha}_{x,\, t}(\mathbb{R}^n\times\mathbb{R})$, then there exists a constant $C>0$ such that
\begin{equation*}
  |u(x,t)-u(y,\tau)|\leq C\left(|x-y|+|t-\tau|^{\frac{1}{2}}\right)^{2\alpha}
\end{equation*}
for any $x,\,y\in\mathbb{R}^n$ and $t,\,\tau\in \mathbb{R}$.
\item[(ii)]
When $\frac{1}{2}<\alpha\leq1$, we say that
$$u(x,t)\in C^{2\alpha,\alpha}_{x,\, t}(\mathbb{R}^n\times\mathbb{R}):=C^{1+(2\alpha-1),\alpha}_{x,\, t}(\mathbb{R}^n\times\mathbb{R}),$$ if $u$ is $\alpha$-H\"{o}lder continuous in $t$ uniformly with respect to $x$ and its gradient $\nabla_xu$ is $(2\alpha-1)$-H\"{o}lder continuous in $x$ uniformly with respect to $t$ and $(\alpha-\frac{1}{2})$-H\"{o}lder continuous in $t$ uniformly with respect to $x$.
\item[(iii)] While for $\alpha>1$, if
$u(x,t)\in C^{2\alpha,\alpha}_{x,\, t}(\mathbb{R}^n\times\mathbb{R}),$
then it means that
$$\partial_tu,\, D^2_xu \in C^{2\alpha-2,\alpha-1}_{x,\, t}(\mathbb{R}^n\times\mathbb{R}).$$
\end{itemize}
In addition, we can analogously define the local parabolic H\"{o}lder space $C^{2\alpha,\alpha}_{x,\, t,\, \rm{loc}}(\mathbb{R}^n\times\mathbb{R})$.

\subsection{Optimal estimate}\label{2.2}

It is well known that applying an integer order derivative to a Schwartz function preserves the regularity.  Now let us discuss whether this property holds for the fractional order operator $(\partial_t-\Delta)_{\rm right}^s$ given in \eqref{master-R}. More precisely, we will proceed to provide the detailed proof for Lemma \ref{lm1} as follows.

\begin{proof}[\bf Proof of Lemma \ref{lm1}\,.]\,
We first show that $(\partial_t-\Delta)_{\rm{right}}^s\varphi(x,t)\in C^{\infty}(\R^n\times \R).$ The change of variables leads to
\begin{eqnarray*}
 (\partial_t-\Delta)_{\rm{right}}^s\varphi(x,t)&=& C_{n,s}\jf_t^{+\infty}\jf_{\R^{n}}\fr{\varphi(x,t)-\varphi(y,\tau)}{(\tau-t)^{\frac{n}{2}+1+s}}e^{-\frac{\abs{x-y}^2}{4(\tau-t)}}\operatorname{d}\!y\operatorname{d}\!\tau
 \\  \nonumber
&=&C_{n,s} \jf_0^{+\infty}\jf_{\R^{n}}\fr{\varphi(x,t)-\varphi(x+\hat{y},t+\hat{\tau})}{\hat{\tau}^{\frac{n}{2}+1+s}}e^{-\frac{\abs{\hat{y}}^2}{4\hat{\tau}}}\operatorname{d}\!\hat{y}\operatorname{d}\!\hat{\tau}.
\end{eqnarray*}
Thereby exchanging the order of differentiation and integration, we derive
\[\pl_x^{\al}\pl_t^\beta\lt((\partial_t-\Delta)_{\rm{right}}^s\varphi(x,t)\rt)=(\partial_t-\Delta)_{\rm{right}}^s\lt(\pl_x^{\al}\pl_t^\beta\varphi(x,t)\rt)\]
for any multi-index $\al=(\al_1,...,\al_n) \mb{~and~} \beta\in\mathbb{N}$.
Hence, it follows from $\varphi(x,t)\in \mathcal{S}(\R^{n}\times \R)$ that
$$(\partial_t-\Delta)_{\rm{right}}^s\varphi(x,t)\in C^{\infty}(\R^n\times \R).$$

In the sequel,  we focus on the estimate of \eqref{1} by  decomposing the whole space $\mathbb{R}^n\tm \mathbb{R}$ into four regions  $\overline{B_1(0)}\tm[-1,1]$, ${B_1^c(0)}\tm[-1,1]$, $\overline{B_1(0)}\tm\lt((-\infty,-1)\cup(1,+\infty)\rt)$ and  ${B_1^c(0)}\tm\lt((-\infty,-1)\cup(1,+\infty)\rt)$.
\\[0.3cm]
\indent
$\mathbf{Case~ 1}$. $(x,t)\in  \overline{B_1(0)}\tm[-1,1].$

Due to $(\partial_t-\Delta)_{\rm{right}}^s\varphi(x,t)\in C^{\infty}(\R^n\times \R)$, it is obvious that
\begin{equation}\label{2}\Abs{(\partial_t-\Delta)_{\rm{right}}^s\varphi(x,t)}\leq C\leq\fr{C}{1+\abs{x}^{n+2+2s}+\abs{t}^{\frac{n}{2}+1+s}} \end{equation}
for $(x,t)\in\overline{B_1(0)}\tm[-1,1]$.
\\[0.3cm]
\indent
$\mathbf{Case~ 2}$. $(x,t)\in{B_1^c(0)}\tm[-1,1].$

In this case, we only need to show that
\begin{equation}\label{3}
\Abs{(\partial_t-\Delta)_{\rm{right}}^s\varphi(x,t)}\leq \fr{C\mathbb{}}{1+\abs{x}^{n+2+2s}}.
\end{equation}
 According to the definition of the nonlocal operator $(\partial_t-\Delta)_{\rm{right}}^s$, we  divide the integral with respect to $y$ into the following two parts by whether the point $x$ is a singularity
  \begin{eqnarray*}
 (\partial_t-\Delta)_{\rm{right}}^s\varphi(x,t)&=& C_{n,s}\jf_t^{+\infty}\jf_{\R^{n}}\fr{\varphi(x,t)-\varphi(y,\tau)}{(\tau-t)^{\frac{n}{2}+1+s}}e^{-\frac{\abs{x-y}^2}{4(\tau-t)}}\operatorname{d}\!y\operatorname{d}\!\tau
 \\  \nonumber
&=& C_{n,s}\jf_t^{+\infty}\lt(\jf_{B_{{\frac{\abs{x}}{2}}}^c(x)}+\jf_{B_{\frac{\abs{x}}{2}}(x)}\rt)\fr{\varphi(x,t)-\varphi(y,\tau)}{(\tau-t)^{\frac{n}{2}+1+s}}e^{-\frac{\abs{x-y}^2}{4(\tau-t)}}
\operatorname{d}\!y\operatorname{d}\!\tau
\\  \nonumber
 &:=&I+II.
 \end{eqnarray*}

With respect to the estimate of $I$, combining Fubini's theorem, $\abs{y-x}\geq\fr{\abs{x}}{2}>\frac{1}{2}$, the fact that
\begin{equation}\label{5}\fr{e^{-\frac{\abs{x-y}^2}{4(\tau-t)}}}{(\tau-t)^{\frac{n}{2}+1+s}}\leq\fr{C}{\Abs{x-y}^{n+2+2s}+(\tau-t)^{\frac{n}{2}+1+s}} \,\, \mbox{for}\,\,\tau>t,\end{equation}
and $\varphi\in \mathcal{S}(\R^n\tm\R)$, we obtain
\begin{eqnarray*}
 \abs{I}&=& \Abs{C_{n,s}\jf_t^{+\infty}\jf_{B_\frac{\abs{x}}{2}^c(x)}\fr{\varphi(x,t)-\varphi(y,\tau)}{(\tau-t)^{\frac{n}{2}+1+s}}e^{-\frac{\abs{x-y}^2}{4(\tau-t)}}\operatorname{d}\!y\operatorname{d}\!\tau
 }\\  \nonumber
&\leq& C_{n,s}\left[\abs{\varphi(x,t)}\jf_{B_\frac{\abs{x}}{2}^c(x)}\lt(\jf_t^{+\infty}\fr{e^{-\frac{\abs{x-y}^2}{4(\tau-t)}}}{(\tau-t)^{\frac{n}{2}+1+s}}d\tau\rt)\operatorname{d}\!y\right.\\
&&\qquad\quad\left.
+\jf_t^{+\infty}\jf_{B_\frac{\abs{x}}{2}^c(x)}\abs{\varphi(y,\tau)}\fr{e^{-\frac{\abs{x-y}^2}{4(\tau-t)}}}{(\tau-t)^{\frac{n}{2}+1+s}}\operatorname{d}\!y\operatorname{d}\!\tau
\right] \\  \nonumber
 &\leq&
  C(n,s)\lt[\abs{\varphi(x,t)}\jf_{B_{\frac{1}{2}}^c(x)}\fr{1}{\abs{x-y}^{n+2s}}\operatorname{d}\!y
+\jf_t^{+\infty}\jf_{B_\frac{\abs{x}}{2}^c(x)}\fr{\abs{\varphi(y,\tau)}}{\abs{x-y}^{n+2+2s}}\operatorname{d}\!y\operatorname{d}\!\tau
\rt] \\  \nonumber
 &\leq&
  C(n,s)\lt[\abs{\varphi(x,t)}
+\fr{1}{\abs{x}^{n+2+2s}}\|\varphi\|_{L^1(\R^n\tm\R)}\rt] \\  \nonumber
 &\overset{}{\leq} &\fr{C(n,s)}{1+\abs{x}^{n+2+2s}}\,.
 \end{eqnarray*}

While for the estimate of $II$, note that the integral  has  singularity at point $(x,t)$ in this situation, so we insert $\varphi(y,t)$ into the integrand $\varphi(x,t)-\varphi(y,\tau)$ to overcome the singularity of the integral. Specifically, we rewrite
\begin{eqnarray*}
II
&=&C_{n,s}\left[\jf_t^{+\infty}\jf_{B_{\frac{\abs{x}}{2}}(x)}\fr{\varphi(x,t)-\varphi(y,t)}{(\tau-t)^{\frac{n}{2}+1+s}}e^{-\frac{\abs{x-y}^2}{4(\tau-t)}}\operatorname{d}\!y\operatorname{d}\!\tau\right.\\
&&\qquad\quad
\left.+\jf_t^{+\infty}\jf_{B_{\frac{\abs{x}}{2}}(x)}\fr{\varphi(y,t)-\varphi(y,\tau)}{(\tau-t)^{\frac{n}{2}+1+s}}e^{-\frac{\abs{x-y}^2}{4(\tau-t)}}\operatorname{d}\!y\operatorname{d}\!\tau\right]
\\
&:=&II_1+II_2.
\end{eqnarray*}

For $II_1$, we apply Taylor expansion for $\varphi(\cdot,t)$ at the point $x$ as follows
\[{\varphi(x,t)-\varphi(y,t)}={\nabla\varphi(x,t)\cdot(x-y)+\frac{1}{2}(x-y)^{T}\nabla^2\varphi(\xi,t)(x-y)},\]
where $\xi$ lies between $x$ and $y$. In views of $y\in B_{\frac{\abs{x}}{2}}(x)$, it implies that the norms of $\xi$, ${x}$ and ${y}$ can be mutually controlled.
Combining the definition of the Cauchy principal value with the fact $\left|\nabla^2\varphi(\xi,t)\right|\leq \frac{C}{1+\left|\xi\right|^{n+4}}$  holds for  $\nabla^2\varphi(\xi,t)\in\mathcal{S}(\R^n\tm\R)$, we derive
\begin{eqnarray*}
\abs{II_1}&=&C_{n,s}\Abs{\jf_{B_{\frac{\abs{x}}{2}}(x)}\lt[\varphi(x,t)-\varphi(y,t)\rt]\jf_t^{+\infty}\fr{e^{-\frac{\abs{x-y}^2}{4(\tau-t)}}}{(\tau-t)^{\frac{n}{2}+1+s}}\operatorname{d}\!\tau \operatorname{d}\!y}
\\
&\leq&C\left|P.V.\jf_{B_{\frac{\abs{x}}{2}}(x)}\fr{\varphi(x,t)-\varphi(y,t)}{\abs{x-y}^{n+2s}}\operatorname{d}\!y\right|
\\
&\leq&C \jf_{B_{\frac{\abs{x}}{2}}(x)}\fr{\Abs{\nabla^2\varphi(\xi,t)}}{\abs{x-y}^{n+2s-2}}\operatorname{d}\!y
\\
&\leq& \fr{C}{1+\abs{x}^{n+4}}\jf_{B_{\frac{\abs{x}}{2}}(x)}\fr{1}{\abs{x-y}^{n+2s-2}}\operatorname{d}\!y
\\
&\leq& \fr{C}{1+\abs{x}^{n+2+2s}}.\end{eqnarray*}

With respect to the estimate of $II_2$, we divide the integral interval of $\tau$ according to the singularity as follows
\begin{eqnarray*}
II_2&=&C_{n,s}\lt[\left(\jf_2^{+\infty}\jf_{B_{\frac{\abs{x}}{2}}(x)}+\jf_t^{2}\jf_{B_{\frac{\abs{x}}{2}}(x)}\right)
\fr{\varphi(y,t)-\varphi(y,\tau)}{(\tau-t)^{\frac{n}{2}+1+s}}e^{-\frac{\abs{x-y}^2}{4(\tau-t)}}\operatorname{d}\!y\operatorname{d}\!\tau
\rt]\\
&:=&II_{2,1}+II_{2,2}.
\end{eqnarray*}
For $II_{2,1}$, we exploit $\varphi\in\mathcal{S}(\R^n\tm\R)$ again to yield that
\begin{eqnarray*}
\Abs{II_{2,1}}&=&C_{n,s}\Abs{\jf_2^{+\infty}\jf_{B_{\frac{\abs{x}}{2}}(x)}\fr{\varphi(y,t)-\varphi(y,\tau)}{(\tau-t)^{\frac{n}{2}+1+s}}e^{-\frac{\abs{x-y}^2}{4(\tau-t)}}\operatorname{d}\!y\operatorname{d}\!\tau
}\\
&\leq&C_{n,s}\jf_2^{+\infty}\jf_{B_{\frac{\abs{x}}{2}}(x)}\lt(\abs{{\varphi(y,t)}}+\abs{\varphi(y,\tau)}\rt)\fr{e^{-\frac{\abs{x-y}^2}{4(\tau-t)}}}{(\tau-t)^{\frac{n}{2}+1+s}}\operatorname{d}\!y\operatorname{d}\!\tau
\\
&\leq&\fr{C}{1+\abs{x}^{n+2+2s}}\jf_2^{+\infty}\jf_{\R^n}\fr{e^{-\frac{\abs{x-y}^2}{4(\tau-t)}}}{(\tau-t)^{\frac{n}{2}+1+s}}\operatorname{d}\!y\operatorname{d}\!\tau
\\
&\leq&\fr{C}{1+\abs{x}^{n+2+2s}}\jf_2^{+\infty}\fr{1}{(\tau-t)^{1+s}}\operatorname{d}\!\tau
\\
&\leq& \fr{C}{1+\abs{x}^{n+2+2s}}.\end{eqnarray*}
While for $II_{2,2}$, we utilize the differential  mean value theorem for $\varphi(y,\cdot)$ near $t$ and the fact that $\pl_t\varphi$ also belongs to $\mathcal{S}(\R^n\tm\R)$,
there exists $\eta\in(t,\tau)$  such that
\[\frac{\abs{\varphi(y,t)-\varphi(y,\tau)}}{\tau-t}=\abs{\pl_t\varphi(y,\eta)}\leq \frac{C}{1+\abs{y}^{n+2+2s}}.\]
Then it follows that
\begin{eqnarray*}
\Abs{II_{2,2}}&=&C_{n,s}\Abs{\jf_t^{2}\jf_{B_{\frac{\abs{x}}{2}}(x)}\fr{\varphi(y,t)-\varphi(y,\tau)}{(\tau-t)^{\frac{n}{2}+1+s}}e^{-\frac{\abs{x-y}^2}{4(\tau-t)}}\operatorname{d}\!y\operatorname{d}\!\tau
}\\
&\leq&\fr{C}{1+\abs{x}^{n+2+2s}}\jf_t^{2}\jf_{\R^n}\fr{e^{-\frac{\abs{x-y}^2}{4(\tau-t)}}}{(\tau-t)^{\frac{n}{2}+s}}\operatorname{d}\!y\operatorname{d}\!\tau
\\
&\leq&\fr{C}{1+\abs{x}^{n+2+2s}}\jf_t^{2}\fr{1}{(\tau-t)^{s}}\operatorname{d}\!\tau
\\
&\leq& \fr{C}{1+\abs{x}^{n+2+2s}}.\end{eqnarray*}
Therefore, we deduce that the assertion \eqref{3} is valid.
\\[0.3cm]
\indent
$\mathbf{Case~ 3}$. $(x,t)\in\overline{B_1(0)}\tm\lt((-\infty,-1)\cup(1,+\infty)\rt).$

In this case, it suffices to prove that
\begin{equation}\label{4}
\Abs{(\partial_t-\Delta)_{\rm{right}}^s\varphi(x,t)}\leq \fr{C\mathbb{}}{1+\abs{t}^{n+2+2s}}.
\end{equation}
By proceeding a similar calculation as in Case 2, so we omit the proof here.
\\[0.3cm]
\indent
$\mathbf{Case~ 4}$. $(x,t)\in{B_1^c(0)}\tm\lt((-\infty,-1)\cup(1,+\infty)\rt).$

At this point, we will only consider the case $(x,t)\in{B_1^c(0)}\tm(-\infty,1)$, as the proof for the other case follows a similar argument.
We first divide the integral domain into four parts as follows
\begin{eqnarray*}
 (\partial_t-\Delta)_{\rm{right}}^s\varphi(x,t)&=& C_{n,s}\jf_t^{+\infty}\jf_{\R^{n}}\fr{\varphi(x,t)-\varphi(y,\tau)}{(\tau-t)^{\frac{n}{2}+1+s}}e^{-\frac{\abs{x-y}^2}{4(\tau-t)}}\operatorname{d}\!y\operatorname{d}\!\tau
 \\  \nonumber
&=& C_{n,s}\lt(\jf_{\frac{t}{2}}^{+\infty}\jf_{B_\frac{\abs{x}}{2}^c(x)}+\jf_t^{\frac{t}{2}}\jf_{B_\frac{\abs{x}}{2}^c(x)}\rt)\fr{\varphi(x,t)-\varphi(y,\tau)}{(\tau-t)^{\frac{n}{2}+1+s}}e^{-\frac{\abs{x-y}^2}{4(\tau-t)}}\operatorname{d}\!y\operatorname{d}\!\tau\\
   &&+\,C_{n,s}\lt(\jf_{\frac{t}{2}}^{+\infty}\jf_{B_\frac{\abs{x}}{2}(x)}+\jf_t^{\frac{t}{2}}\jf_{B_\frac{\abs{x}}{2}^c(x)}\rt)\fr{\varphi(x,t)-\varphi(y,\tau)}{(\tau-t)^{\frac{n}{2}+1+s}}e^{-\frac{\abs{x-y}^2}{4(\tau-t)}}\operatorname{d}\!y\operatorname{d}\!\tau
\\  \nonumber
 &:=&I+II+III+IV.
 \end{eqnarray*}

Regarding the estimate of $I$, we utilize \eqref{5},  $\varphi\in \mathcal{S}(\R^n\tm\R)$, and $|x|,\,|t|>1$ in this case.  By incorporating these conditions, we derive
\begin{eqnarray*}
 \Abs{I}&=&\Abs{C_{n,s}\jf_{\frac{t}{2}}^{+\infty}\jf_{B_\frac{\abs{x}}{2}^c(x)}\fr{\varphi(x,t)-\varphi(y,\tau)}{(\tau-t)^{\frac{n}{2}+1+s}}e^{-\frac{\abs{x-y}^2}{4(\tau-t)}}\operatorname{d}\!y\operatorname{d}\!\tau
}\\ &\leq&C_{n,s}\left[\abs{\varphi(x,t)}\jf_{\frac{t}{2}}^{+\infty}\jf_{B_\frac{\abs{x}}{2}^c(x)}\fr{e^{-\frac{\abs{x-y}^2}{4(\tau-t)}}}{(\tau-t)^{\frac{n}{2}+1+s}}\operatorname{d}\!y\operatorname{d}\!\tau\right.\\
&&\qquad\quad\left.
+\jf_{\frac{t}{2}}^{+\infty}\jf_{B_\frac{\abs{x}}{2}^c(x)}\abs{\varphi(y,\tau)}\fr{e^{-\frac{\abs{x-y}^2}{4(\tau-t)}}}{(\tau-t)^{\frac{n}{2}+1+s}}\operatorname{d}\!y\operatorname{d}\!\tau
\right]\\
&\overset{}{\leq}&C\left[\abs{\varphi(x,t)}\jf_{\frac{t}{2}}^{+\infty}\fr{1}{(\tau-t)^{1+s}}\operatorname{d}\!\tau\right.\\
&&\qquad\quad\left.+\jf_{\frac{t}{2}}^{+\infty}\jf_{B_\frac{\abs{x}}{2}^c(x)}\abs{\varphi(y,\tau)}\fr{1}{\Abs{x-y}^{n+2+2s}+(\tau-t)^{\frac{n}{2}+1+s}}\operatorname{d}\!y\operatorname{d}\!\tau
\right]\\
&\leq&C\lt[\abs{\varphi(x,t)}\fr{1}{\abs{t}^s}
+\fr{1}{\Abs{x}^{n+2+2s}+\abs{t}^{\frac{n}{2}+1+s}}\|\varphi\|_{L^1(\R^n\tm\R)}
\rt]\\
 &\overset{}{\leq}&\fr{C}{1+\Abs{x}^{n+2+2s}+\abs{t}^{\frac{n}{2}+1+s}}.
\end{eqnarray*}

To estimate $II$, combining \eqref{5} with $$\abs{\varphi(y,\tau)}\leq\fr{C}{\lt(1+\Abs{y}^{k}\rt)\lt(1+\abs{\tau}^{\frac{n}{2}+2+s}\rt)} \mb{~for~some~} k>n$$
due to $\varphi\in\mathcal{S}$, we calculate
\begin{eqnarray*}
 \Abs{II}
&=& C_{n,s}\Abs{\jf_t^{\frac{t}{2}}\jf_{B_\frac{\abs{x}}{2}^c(x)}\fr{\varphi(x,t)-\varphi(y,\tau)}{(\tau-t)^{\frac{n}{2}+1+s}}e^{-\frac{\abs{x-y}^2}{4(\tau-t)}}\operatorname{d}\!y\operatorname{d}\!\tau}\\
   &\leq&C\left[\abs{\varphi(x,t)}\jf_{B_\frac{\abs{x}}{2}^c(x)}\jf_t^{+\infty}\fr{e^{-\frac{\abs{x-y}^2}{4(\tau-t)}}}{(\tau-t)^{\frac{n}{2}+1+s}}\operatorname{d}\!\tau \operatorname{d}\!y\right.\\
   &&\qquad\quad\left.+\jf_t^{\frac{t}{2}}\jf_{B_\frac{\abs{x}}{2}^c(x)}\abs{\varphi(y,\tau)}\fr{1}{\Abs{x-y}^{n+2+2s}}\operatorname{d}\!y\operatorname{d}\!\tau
\right]\\
 &\overset{}{\leq}&C\lt[\abs{\varphi(x,t)}\jf_{B_\frac{1}{2}^c(x)}\fr{1}{\abs{x-y}^{n+2s}}\operatorname{d}\!y
+\fr{1}{\Abs{x}^{n+2+2s}}\jf_t^{\frac{t}{2}}\jf_{B_\frac{\abs{x}}{2}^c(x)}\abs{\varphi(y,\tau)}\operatorname{d}\!y\operatorname{d}\!\tau
\rt]\\
&\overset{}{\leq}&C\lt[\fr{1}{1+\abs{x}^{n+2+2s}+\abs{t}^{\frac{n}{2}+1+s}}
+\fr{1}{\Abs{x}^{n+2+2s}}\jf_{\R^n}\fr{1}{1+\abs{y}^k}dy\jf_t^{\frac{t}{2}}\fr{1}{1+|\tau|^{\frac{n}{2}+2+s}}\operatorname{d}\!\tau
\rt]\\
&\leq&C\lt[\fr{1}{1+\abs{x}^{n+2+2s}+\abs{t}^{\frac{n}{2}+1+s}}
+\fr{1}{\Abs{x}^{n+2+2s}(1+\abs{t}^{\frac{n}{2}+1+s})}
\rt]\\
&\overset{}{\leq}&\fr{C}{1+\abs{x}^{n+2+2s}+\abs{t}^{\frac{n}{2}+1+s}}\,.
 \end{eqnarray*}

For $III$,  in analogy with the estimate of $II$, we apply $|t|>1$ and
 $$\abs{\varphi(y,\tau)}\leq\fr{C}{\lt(1+\Abs{y}^{2n+2+2s}\rt)\lt(1+\abs{\tau}^{m}\rt)} \mb{~for~some~} m>1$$
to arrive at
 \begin{eqnarray*}
 \Abs{III}
&=& C_{n,s}\Abs{\jf_{\frac{t}{2}}^{+\infty}\jf_{B_\frac{\abs{x}}{2}(x)}\fr{\varphi(x,t)-\varphi(y,\tau)}{(\tau-t)^{\frac{n}{2}+1+s}}e^{-\frac{\abs{x-y}^2}{4(\tau-t)}}\operatorname{d}\!y\operatorname{d}\!\tau}\\
   &\overset{}{\leq}&C\lt[\abs{\varphi(x,t)}\jf_{\frac{t}{2}}^{+\infty}\jf_{B_\frac{\abs{x}}{2}(x)}\fr{e^{-\frac{\abs{x-y}^2}{4(\tau-t)}}}{(\tau-t)^{\frac{n}{2}+1+s}}\operatorname{d}\!y\operatorname{d}\!\tau
+\jf_{\frac{t}{2}}^{+\infty}\jf_{B_\frac{\abs{x}}{2}(x)}\fr{\abs{\varphi(y,\tau)}}{(\tau-t)^{\frac{n}{2}+1+s}}\operatorname{d}\!y\operatorname{d}\!\tau
\rt]\\
 &{\leq}&C\lt[\abs{\varphi(x,t)}\jf_{\frac{t}{2}}^{+\infty}\fr{1}{(\tau-t)^{1+s}}\operatorname{d}\!\tau
+\fr{1}{\Abs{t}^{\frac{n}{2}+1+s}}\jf_{\frac{t}{2}}^{+\infty}\jf_{B_\frac{\abs{x}}{2}(x)}\abs{\varphi(y,\tau)}\operatorname{d}\!y\operatorname{d}\!\tau
\rt]\\
&\overset{}{\leq}&C\lt[\fr{1}{1+\abs{x}^{n+2+2s}+\abs{t}^{\frac{n}{2}+1+s}}\fr{1}{\abs{t}^{s}}
+\fr{1}{\Abs{t}^{\frac{n}{2}+1+s}}\jf_{0}^{+\infty}\fr{1}{1+\tau^m}d\tau\jf_{B_\frac{\abs{x}}{2}(x)}\fr{1}{1+\abs{y}^k}dy
\rt]\\
&\underset{}{\overset{}{\leq}}&C\lt[\fr{1}{1+\abs{x}^{n+2+2s}+\abs{t}^{\frac{n}{2}+1+s}}
+\fr{1}{\Abs{t}^{\frac{n}{2}+1+s}(1+\abs{x}^{n+2+2s})}
\rt]\\
&\overset{}{\leq}&\fr{C}{1+\abs{x}^{n+2+2s}+\abs{t}^{\frac{n}{2}+1+s}}.
 \end{eqnarray*}

 Next, we estimate $IV$, which has singularity  at the point $(x,t)$. Dividing the integral into two parts as in Case 2
 \begin{eqnarray*}
 IV
& =&C_{n,s}\jf_t^{\frac{t}{2}}\jf_{B_\frac{\abs{x}}{2}(x)}\fr{\varphi(x,t)-\varphi(y,t)}{(\tau-t)^{\frac{n}{2}+1+s}}e^{-\frac{\abs{x-y}^2}{4(\tau-t)}}dyd\tau\\
&&+\, C_{n,s}\jf_t^{\frac{t}{2}}\jf_{B_\frac{\abs{x}}{2}(x)}\fr{\varphi(y,t)-\varphi(y,\tau)}{(\tau-t)^{\frac{n}{2}+1+s}}e^{-\frac{\abs{x-y}^2}{4(\tau-t)}}\operatorname{d}\!y\operatorname{d}\!\tau\\
\\
&:=&IV_1+IV_2.
 \end{eqnarray*}
 Similar to the estimate of $II_1$ in Case 2, we can deduce that
 \begin{eqnarray*}
 \Abs{IV_1}
  &\leq& \fr{C}{1+\abs{x}^{n+2+2s}+\abs{t}^{\frac{n}{2}+1+s}}.
\end{eqnarray*}
Moreover, following  a similar calculation of $II_{2,2}$ as in Case 2, we utilize the differential mean value theorem for $\varphi(y,\cdot)$ near $t$ and the fact $\pl_t\varphi\in\mathcal{S}(\R^n\tm\R)$ to obtain
\[\frac{\abs{\varphi(y,t)-\varphi(y,\tau)}}{\tau-t}
\leq
\frac{C}{1+\abs{y}^{n+2+2s}+\abs{t}^{\frac{n}{2}+2}},\] which yields that
\begin{eqnarray*}
\Abs{IV_{2}}
&\leq&\fr{C}{1+\abs{x}^{n+2+2s}+\abs{t}^{\frac{n}{2}+2}}\jf_t^{\frac{t}{2}}\jf_{\R^n}\fr{e^{-\frac{\abs{x-y}^2}{4(\tau-t)}}}{(\tau-t)^{\frac{n}{2}+s}}\operatorname{d}\!y\operatorname{d}\!\tau
\\
&\leq&\fr{C}{1+\abs{x}^{n+2+2s}+\abs{t}^{\frac{n}{2}+2}}\jf_t^{\frac{t}{2}}\fr{1}{(\tau-t)^{s}}\operatorname{d}\!\tau
\\ &\leq& \fr{C}{1+\abs{x}^{n+2+2s}+\abs{t}^{\frac{n}{2}+1+s}}.
\end{eqnarray*}
Hence, we conclude that the estimate \eqref{1} holds in this case.

In conclusion, together with the above four cases, we establish the validity of estimate \eqref{1}. Thus, the proof of Lemma \ref{lm1} is complete.
\end{proof}

Lemma \ref{lm1} clearly demonstrates that $(\partial_t-\Delta)_{\rm{right}}^s\varphi(x,t)$ does not belong to $\mathcal{S}(\R^{n}\times \R)$ when $\varphi(x,t)\in \mathcal{S}(\R^{n}\times \R)$. Despite it loss of rapid decay at any order, $(\partial_t-\Delta)_{\rm{right}}^s\varphi$ still retains the smoothness and exhibits decay with a specific exponent.
To illustrate the optimality of the decay exponent in \eqref{1}, let us consider the following counterexample.

\begin{counterexample}
Let $\eta\in  C_0^\infty(\R^n\tm\R)$ satisfying $-1\leq\eta\leq0$ and
 \begin{equation*}\eta(x, t)  = \begin{cases} -1, & (x,t)\in \overline{B_1(0)}\times[-1,1], \\ 0, & (x,t)\notin B_2(0)\times(-4,4). \end{cases}\end{equation*}
 If $\abs{x}^2\sim\abs{t}$ with sufficiently negative $t$, then we derive
  \begin{eqnarray*}
\Abs{(\partial_t-\Delta)_{\rm{right}}^s\eta(x,t)}&=&C_{n,s}\jf_t^{+\infty}\jf_{\R^n}\fr{-\eta(y,\tau)}{(\tau-t)^{\frac{n}{2}+1+s}}
e^{-\frac{\abs{x-y}^2}{4(\tau-t)}}\operatorname{d}\!t\operatorname{d}\!\tau
\\
&\geq&C_{n,s}\jf_{-1}^1\jf_{B_1(0)}\fr{1}{(\tau-t)^{\frac{n}{2}+1+s}}e^{-\frac{\abs{x-y}^2}{4(\tau-t)}}\operatorname{d}\!t\operatorname{d}\!\tau\\
&\geq&\fr{C}{|t|^{\frac{n}{2}+1+s}}e^{-\frac{\abs{x}^2}{4|t|}}\\
&\geq&\fr{C_0}{1+\abs{x}^{n+2+2s}+|t|^{\frac{n}{2}+1+s}},
\end{eqnarray*}
where $C_0$ is a positive constant that depends only on $n$ and $s$.
This shows that  estimate \eqref{1} cannot be improved on the set $\left\{(x,t)\in\R^n\tm\R\mid  \abs{x}^2\sim\abs{t}\right\}$. Indeed, if not, suppose that there exist some positive constants $C_1$ and $\delta_1$ such that
\[\Abs{(\partial_t-\Delta)_{\rm{right}}^s\eta(x,t)}\leq\fr{C_1}{1+\abs{x}^{n+2+2s}+\abs{t}^{\frac{n}{2}+1+s+\delta_1}}
~\mb{in}~\left\{(x,t)\in\R^n\tm\R\mid  \abs{x}^2\sim\abs{t}\right\}.\]
Then  for sufficiently negative $t$  as aforementioned,  we have
\[\fr{C_0}{1+\abs{x}^{n+2+2s}+|t|^{\frac{n}{2}+1+s}}\leq\Abs{(\partial_t-\Delta)_{\rm{right}}^s\eta(x,t)}\leq\fr{C_1}{1+\abs{x}^{n+2+2s}
+\abs{t}^{\frac{n}{2}+1+s+\delta_1}}.\]
It follows from $(x,t)\in \left\{(x,t)\in\R^n\tm\R \mid  \abs{x}^2\sim\abs{t}\right\}$ that
\[\fr{C_1}{C_0}\geq\fr{1+\abs{x}^{n+2+2s}+|t|^{\frac{n}{2}+1+s+\delta_1}}{1+\abs{x}^{n+2+2s}+|t|^{\frac{n}{2}+1+s}}\geq C(1+\abs{t}^{\delta_1}),\]
which is a contradiction as $t\to-\infty$. In analogy with above estimates, we can also derive a contradiction if there exist some positive constants $C_2$ and $\delta_2$ such that
\[\Abs{(\partial_t-\Delta)_{\rm{right}}^s\eta(x,t)}\leq\fr{C_2}{1+\abs{x}^{n+2+2s+\delta_2}+\abs{t}^{\frac{n}{2}+1+s}}
~\mb{in}~\left\{(x,t)\in\R^n\tm\R\mid  \abs{x}^2\sim\abs{t}\right\}.\]
Therefore, this counterexample demonstrates that the decay exponent in \eqref{1} is optimal.
\end{counterexample}

\subsection{Fractional derivative and Fourier transform of tempered distributions}\label{2.3}
In this subsection, we first define $(\partial_t-\Delta)^s$ in the sense of distributions.
Note that if $u$ and $\varphi$ are both Schwartz functions, it follows from Parseval's relation that
\begin{eqnarray}\label{derivative}
    &&\int_{-\infty}^{+\infty} \int_{\mathbb{R}^n} (\partial_t-\Delta)^su(x,t)\overline{\varphi(x,t)}\operatorname{d}\!x\operatorname{d}\!t\nonumber\\
   &=& \int_{-\infty}^{+\infty} \int_{\mathbb{R}^n} \mathcal{F}\left((\partial_t-\Delta)^su\right)(\xi,\rho)\overline{\mathcal{F}\varphi(\xi,\rho)}\operatorname{d}\!\xi\operatorname{d}\!\rho \nonumber \\
   &=&  \int_{-\infty}^{+\infty} \int_{\mathbb{R}^n} (i\rho+|\xi|^2)^s \mathcal{F}u(\xi,\rho)\overline{\mathcal{F}\varphi(\xi,\rho)}\operatorname{d}\!\xi\operatorname{d}\!\rho\nonumber\\
    &=&  \int_{-\infty}^{+\infty} \int_{\mathbb{R}^n} \mathcal{F}u(\xi,\rho)\overline{(-i\rho+|\xi|^2)^s \mathcal{F}\varphi(\xi,\rho)}\operatorname{d}\!\xi\operatorname{d}\!\rho\nonumber\\
    &=&  \int_{-\infty}^{+\infty} \int_{\mathbb{R}^n} \mathcal{F}u(\xi,\rho)\overline{\mathcal{F}\left((\partial_t-\Delta)_{\rm right}^s\varphi\right)(\xi,\rho)}\operatorname{d}\!\xi\operatorname{d}\!\rho\nonumber\\
       &=&  \int_{-\infty}^{+\infty} \int_{\mathbb{R}^n} u(x,t)\overline{(\partial_t-\Delta)_{\rm right}^s\varphi(x,t)}\operatorname{d}\!x\operatorname{d}\!t.
\end{eqnarray}
To justify the above calculation, it suffices to show that
\begin{equation}\label{fourier-R}
  \mathcal{F}\left((\partial_t-\Delta)_{\rm right}^s \varphi\right)(\xi,\rho)=(-i\rho+|\xi|^2)^s\mathcal{F}\varphi(\xi,\rho)
\end{equation}
for $\varphi\in \mathcal{S}(\mathbb{R}^n\times\mathbb{R})$.

The fully fractional heat operator $(\partial_t-\Delta)_{\rm right}^s $ defined in \eqref{master-R} can be equivalently expressed as follows
\begin{equation*}
(\partial_t-\Delta)_{\rm right}^s \varphi(x,t)
=\frac{ C_{n,s}}{2}\int_{0}^{+\infty}\int_{\mathbb{R}^n}
  \frac{2\varphi(x,t)-\varphi(x+y,t+\tau)-\varphi(x-y,t+\tau)}
  {\tau^{\frac{n}{2}+1+s}}e^{-\frac{|y|^2}{4\tau}}\operatorname{d}\!y\operatorname{d}\!\tau.
\end{equation*}
This representation does not require taking the Cauchy principal value with respect to $y$, enabling us to utilize the Fubini's theorem to derive
\begin{eqnarray*}
   &&  \mathcal{F}\left((\partial_t-\Delta)_{\rm right}^s \varphi\right)(\xi,\rho) \\
   &=& \frac{ C_{n,s}}{2}\int_{0}^{+\infty}\int_{\mathbb{R}^n}
  \frac{2\mathcal{F}(\varphi(x,t))-\mathcal{F}(\varphi(x+y,t+\tau))-\mathcal{F}(\varphi(x-y,t+\tau))}
  {\tau^{\frac{n}{2}+1+s}}e^{-\frac{|y|^2}{4\tau}}\operatorname{d}\!y\operatorname{d}\!\tau \\
   &=&   C_{n,s}\mathcal{F}\varphi(\xi,\rho)\int_{0}^{+\infty}\int_{\mathbb{R}^n}
  \frac{1-e^{i(y\cdot\xi+\tau\rho)}}
  {\tau^{\frac{n}{2}+1+s}}e^{-\frac{|y|^2}{4\tau}}\operatorname{d}\!y\operatorname{d}\!\tau \\
  &=&  \frac{1}{\Gamma(-s)} \mathcal{F}\varphi(\xi,\rho) \int_{0}^{+\infty}\frac{e^{-\tau(|\xi|^2-i\rho)}-1}{\tau^{1+s}}\operatorname{d}\!\tau\\
  &=&(-i\rho+|\xi|^2)^s\mathcal{F}\varphi(\xi,\rho),
\end{eqnarray*}
where the Cauchy integral theorem and  the definition of Gamma function are applied in the last line. Hence, we verify that \eqref{fourier-R} is valid.

From \eqref{fourier-R}, we obverse that $\mathcal{F}\left((\partial_t-\Delta)_{\rm{right}}^s\varphi\right)\notin\mathcal{S}$ since $(-i\rho+\abs{\xi}^2)^s$ exhibits a singularity at the origin in its Fourier transform, similar to the Fractional laplacian. This provides further confirmation that the result we obtained in Lemma \ref{lm1} regarding $(\partial_t-\Delta)_{\rm{right}}^s\varphi\notin\mathcal{S}$ is correct. Meanwhile, Lemma \ref{lm1} shows that  $(\partial_t-\Delta)_{\rm{right}}^s\varphi\in \widetilde{\mathcal{S}}$ when $\varphi\in\mathcal{S}$, where $\widetilde{\mathcal{S}}\supset\mathcal{S}$ denoted by
   \[\widetilde{\mathcal{S}}(\R^n\tm\R):=\left\{u\in C^{\infty}(\R^n\tm\R)\mid\rho_{\alpha,\beta}(u)<+\infty,\ \forall \mb{~multi-index~} \alpha \mb{~and~} \beta\in \mathbb{N}\right\},\]
where the topology in $\widetilde{\mathcal{S}}$
is given by the
family of seminorms$$\rho_{\alpha,\beta}(u):=\sup\limits_{(x,t)\in \R^n\tm\R}(1+\abs{x}^{n+2+2s}+\abs{t}^{\frac{n}{2}+1+s})|\partial_x^\alpha \partial_t^\beta u(x,t)|.$$
 We take $\widetilde{\mathcal{S}}'$ to be the dual of $\widetilde{\mathcal{S}}$.  Now by duality as in \eqref{derivative}, we can extend the definition of $(\partial_t-\Delta)^s$ in $\mathcal{S}$ to the
  space $\widetilde{\mathcal{S}}'$ of tempered distributions. In general, we define the operator $(\partial_t-\Delta)^s$ in a natural space we are going to use, that is,
  \begin{eqnarray*}
    \mathcal{L}^{2s,s}(\mathbb{R}^n\times\mathbb{R}) &:=& L^1_{\rm loc} (\mathbb{R}^n\times\mathbb{R})\cap\widetilde{\mathcal{S}}' \\
    &=&  \left\{u:\R^n\tm\R\to\R\mid \int_{-\infty}^{+\infty} \int_{\mathbb{R}^n} \frac{|u(x,t)|}{1+|x|^{n+2+2s}+|t|^{\frac{n}{2}+1+s}}\operatorname{d}\!x\operatorname{d}\!t<\infty\right\}.
  \end{eqnarray*}
It is worth  emphasizing that the space $\mathcal{L}^{2s,s}(\mathbb{R}^n\times\mathbb{R})$ is clearly a subset of the space $\mathcal{L}(\mathbb{R}^n\times\mathbb{R})$ mentioned in Section \ref{sec1}\,.

Now we define the left fractional derivative of a distribution as follows.
\begin{definition}\label{dlfd}
Let $u\in \mathcal{L}^{2s,s}(\mathbb{R}^n\times\mathbb{R})$, we define $(\partial_t-\Delta)^su(x,t)$ by
\begin{equation*}
  \int_{-\infty}^{+\infty} \int_{\mathbb{R}^n} (\partial_t-\Delta)^su(x,t)\overline{\varphi(x,t)}\operatorname{d}\!x\operatorname{d}\!t=\int_{-\infty}^{+\infty} \int_{\mathbb{R}^n} u(x,t)\overline{(\partial_t-\Delta)_{\rm right}^s\varphi(x,t)}\operatorname{d}\!x\operatorname{d}\!t
\end{equation*}
for any $\varphi\in \mathcal{S}(\mathbb{R}^n\times\mathbb{R})$.
\end{definition}

At the end of this subsection, we recall the definition of the Fourier transform of a tempered distribution $u\in \mathcal{L}^{2s,s}(\mathbb{R}^n\times\mathbb{R})\subset\mathcal{S}'$, as introduced in  \cite{L}.
\begin{definition}\label{dft}
Let $u\in \mathcal{L}^{2s,s}(\mathbb{R}^n\times\mathbb{R})$, we define the Fourier transform $\mathcal{F}u$ of a tempered distribution by
\begin{equation*}
  \langle\mathcal{F}u,\varphi\rangle:=\langle u,\mathcal{F}^{-1}\varphi\rangle\end{equation*}
for all $\varphi\in \mathcal{S}(\mathbb{R}^n\times\mathbb{R})$, where \[\langle u,\mathcal{F}^{-1}\varphi\rangle:=\int_{-\infty}^{+\infty} \int_{\mathbb{R}^n}u(x,t)\overline{\mathcal{F}^{-1}\varphi(x,t)}\operatorname{d}\!x\operatorname{d}\!t.\]
\end{definition}

\section{The proof of main results}\label{sec3}

In this section, we complete the proof of Liouville theorem for the homogeneous master equation (Theorem \ref{Liouville}). As a byproduct, we establish an equivalence between the master equation \eqref{pdeq}
and the corresponding integral equation \eqref{ieq} (Theorem \ref{equivalence}). The significance of this equivalence lies in the fact that it allows us to transform the study of qualitative properties of solutions to the fully fractional heat equations, which may be challenging, into the study of its equivalent integral equations.

\subsection{The proof of Liouville theorem}\label{3.1}

\begin{proof}
[\bf Proof of Theorem \ref{Liouville}\,.] \,
A combination of the definition of distributional solutions \eqref{dissol} and the fact \eqref{fourier-R} yields that
\begin{eqnarray}\label{liou1}
   0&=&\int_{-\infty}^{+\infty} \int_{\mathbb{R}^n} u(x,t)\overline{(\partial_t-\Delta)_{\rm right}^s\varphi(x,t)}\operatorname{d}\!x\operatorname{d}\!t\nonumber \\
  &=& \int_{-\infty}^{+\infty} \int_{\mathbb{R}^n}u(x,t)\overline{\mathcal{F}^{-1}\left((-i\rho+|\xi|^2)^s \mathcal{F}\varphi(\xi,\rho)\right)(x,t)}\operatorname{d}\!x\operatorname{d}\!t
\end{eqnarray}
for any $\varphi\in \mathcal{S}(\mathbb{R}^n\times\mathbb{R})$.

Since $u\in\mathcal{L}^{2s,s}(\mathbb{R}^n\times\mathbb{R})$ is a tempered distribution, then $u$ admits the Fourier transform in the sense of distributions. To proceed, we claim that
\begin{equation}\label{liou2}
  \langle\mathcal{F}u,\phi\rangle=0\,\,\mbox{for any}\,\, \phi\in C_0^\infty(\mathbb{R}^n\times\mathbb{R}\setminus \{({\bf 0}, 0)\}).
\end{equation}
Let $\phi\in C_0^\infty(\mathbb{R}^n\times\mathbb{R}\setminus \{({\bf 0}, 0)\})$, then the function $\frac{\phi(\xi,\rho)}{(-i\rho+|\xi|^2)^s }$ also belongs to $C_0^\infty(\mathbb{R}^n\times\mathbb{R}\setminus \{({\bf 0}, 0)\})$. There must exists a function $\varphi\in \mathcal{S}(\mathbb{R}^n\times\mathbb{R})$ such that
\begin{equation*}
  \mathcal{F}\varphi(\xi,\rho)=\frac{\phi(\xi,\rho)}{(-i\rho+|\xi|^2)^s}.
\end{equation*}
It follows from Definition \ref{dft} and \eqref{liou1} that
\begin{eqnarray*}
 \langle\mathcal{F}u,\phi\rangle&=& \langle\mathcal{F}u,(-i\rho+|\xi|^2)^s\mathcal{F}\varphi\rangle \\
  &=& \int_{-\infty}^{+\infty} \int_{\mathbb{R}^n}u(x,t)\overline{\mathcal{F}^{-1}\left((-i\rho+|\xi|^2)^s \mathcal{F}\varphi(\xi,\rho)\right)(x,t)}\operatorname{d}\!x\operatorname{d}\!t=0.
\end{eqnarray*}
Hence, the assertion \eqref{liou2} holds, which implies that the support of
$\mathcal{F}u$ is a point, the origin. From this, we conclude that $u(x,t)$ is a polynomial of $x$ and $t$. While the condition $u\in\mathcal{L}^{2s,s}(\mathbb{R}^n\times\mathbb{R})$ with $0<s <1$
implies that this polynomial must be of degree one in the form:
$$u(x,t)=C_0+\sum_{i=1}^{n}C_ix_i,$$
where $C_i$ with $i=0,\,1,\,...,\,n$ are constants.
In the case $0 < s \leq \frac{1}{2}$, we directly deduce that $C_i = 0$ for $i=1,\,2,\,...,\,n$ by $u\in\mathcal{L}^{2s,s}(\mathbb{R}^n\times\mathbb{R})$. While in the case $\frac{1}{2} < s < 1$, the asymptotic condition \eqref{AA} of $u$ also implies that $C_i = 0$ for $i=1,\,2,\,...,\,n$.
Hence, we must have
$$u(x,t)\equiv C.$$
This completes the proof of Theorem \ref{Liouville}\,.
\end{proof}

\subsection{The proof of equivalence between pseudo-differential equation and integral equation}\label{3.2}

\begin{proof}
[\bf Proof of Theorem \ref{equivalence}\,.] \,
We start by assuming that $u$ is a solution to the pseudo-differential equation \eqref{pdeq}, where the nonnegative continuous non-homogeneous term $f$ fulfills
\begin{equation}\label{fcondition}
  f(x,t,\zeta)\geq C_2>0\,\, \mbox{for} \,\, \zeta\geq C_1>0,\,\, \mbox{uniformly in} \,\, (x,t)\in \mathbb{R}^n\times \mathbb{R}.
\end{equation}
We aim to prove that $u$ is also a solution to the integral equation \eqref{ieq}.

In order to demonstrate this, we let the smooth cut-off functions
$\phi(x)\in [0,1]$ and $\eta(t)\in[0,1]$ satisfy
\begin{subequations}
\begin{align*}
\phi(x)=
\begin{aligned}
\left\lbrace
\begin{aligned}
1, & \ \ x\in B_{\frac{1}{2}}(0) \,, \\
0, & \ \ x\not\in B_{1}(0)\,,
\end{aligned}
\right.
\end{aligned}
\quad\quad
\eta(t)=
\begin{aligned}
\left\lbrace
\begin{aligned}
1, & \ \ t\in(-\frac{1}{2},\frac{1}{2}) \,, \\
0, & \ \ t\not\in (-1,1)\,.
\end{aligned}
\right.
\end{aligned}
\end{align*}
\end{subequations}
Meanwhile, both $\phi$ and $\eta$ possess radial decreasing behavior with respect to the origin. We denote
\begin{equation*}
  \phi_R(x)=\phi\big(\frac{x}{R}\big)\,\, \mbox{and}\,\, \eta_R(t)=\phi\big(\frac{t}{R^2}\big),
\end{equation*}
then based on the conditions of $\phi$ and $\eta$, it is easy to observe that both $\phi_R$ and $\eta_R$ are monotonically increasing with respect to $R$ for any fixed point $(x,t)$.
Now we set the truncation of $f$ as
\begin{equation*}
  f_R(x,t,u):=f(x,t,u(x,t))\phi_R(x)\eta_R(t).
\end{equation*}
According to the definition of the cut-off functions $\phi_R$ and $\eta_R$, as well as the non-negativity of $f$, it follows that $f_R$ exhibits a monotonically increasing behavior with respect to $R$ for any fixed point $(x,t)$. More precisely, $f_R$ increases towards $f$ as $R\rightarrow\infty$. Furthermore, the
support set of $f_R$ belongs to $B_R(0)\times(-R^2,R^2)$, then by virtue of the continuity of $f$, we deduce that
\begin{equation}\label{supf_R}
  M_R:=\sup_{\mathbb{R}^n\times\mathbb{R}}f_R(x,t,u(x,t))<\infty.
\end{equation}
We define
\begin{equation*}
  v_R(x,t):=\int_{-\infty}^{+\infty} \int_{\mathbb{R}^n} G(x-y,t-\tau)f_R(y,\tau,u(y,\tau))\operatorname{d}\!y\operatorname{d}\!\tau,
\end{equation*}
then \cite{ST} reveals that $v_R$ is a solution of
\begin{equation*}
  (\partial_t-\Delta)^s v_R(x,t) =f_R(x,t,u(x,t)),\,\,(x, t)\in \mathbb{R}^n\times \mathbb{R}.
\end{equation*}
Let $$w_R(x,t):=u(x,t)-v_R(x,t),$$
it is obvious that $w_R$ satisfies
\begin{equation}\label{w_R}
  (\partial_t-\Delta)^s w_R(x,t) =f(x,t,u(x,t))-f_R(x,t,u(x,t))\geq 0,\,\,\mbox{in}\,\, \mathbb{R}^n\times \mathbb{R}.
\end{equation}
In the following, our objective is to establish that
\begin{equation}\label{mp}
  w_R(x,t)\geq 0\,\,\mbox{in}\,\,\mathbb{R}^n\times \mathbb{R}.
\end{equation}
The proof goes by contradiction, if \eqref{mp} is violated, then there exits some point $(x,t)\in\mathbb{R}^n\times \mathbb{R}$ such that $w_R(x,t)<0$.
To derive a contradiction at the negative infimum of $w_R$, we need to show that the negative infimum of $w_R$ is necessarily attained in the interior.
For this purpose, it suffices to claim that the validity of
\begin{equation}\label{w_Rd1}
\varliminf_{|x|\rightarrow\infty} w_R(x,t)\geq 0 \,\,\mbox{uniformly for}\,\,t\in\mathbb{R},
\end{equation}
and
\begin{equation}\label{w_Rd2}
\varliminf_{|t|\rightarrow\infty} w_R(x,t)\geq 0 \,\,\mbox{uniformly for}\,\,x\in\mathbb{R}^n.
\end{equation}
Indeed, due to the positivity of $u$, if we demonstrate that
\begin{equation}\label{v_Rd1}
\lim_{|x|\rightarrow\infty} v_R(x,t)= 0 \,\,\mbox{uniformly for}\,\,t\in\mathbb{R},
\end{equation}
and
\begin{equation}\label{v_Rd2}
\lim_{|t|\rightarrow\infty} v_R(x,t)= 0 \,\,\mbox{uniformly for}\,\,x\in\mathbb{R}^n,
\end{equation}
then the aforementioned assertions are evidently true.
We start by establishing the validity of \eqref{v_Rd1}. By virtue of \eqref{supf_R}, and combining Fubini's theorem with Lebesgue's dominated convergence theorem, we arrive at
\begin{eqnarray*}
 \lim_{|x|\rightarrow\infty} v_R(x,t)
  &=& \lim_{|x|\rightarrow\infty} \int_{-\infty}^{+\infty} \int_{\mathbb{R}^n} G(x-y,t-\tau)f_R(y,\tau,u(y,\tau))\operatorname{d}\!y\operatorname{d}\!\tau \\
  &=& A_{n,s}\lim_{|x|\rightarrow\infty}\int_{-R^2}^{t} \int_{B_R(0)} \frac{e^{-\frac{|x-y|^2}{4(t-\tau)}}}{(t-\tau)^{\frac{n}{2}+1-s}}
  f_R(y,\tau,u(y,\tau))\operatorname{d}\!y\operatorname{d}\!\tau\\
   &\leq& A_{n,s}M_R\lim_{|x|\rightarrow\infty} \int_{B_R(0)}\int_{-\infty}^{t} \frac{e^{-\frac{|x-y|^2}{4(t-\tau)}}}{(t-\tau)^{\frac{n}{2}+1-s}}\operatorname{d}\!\tau
  \operatorname{d}\!y\\
  &=& C(n,s)M_R\lim_{|x|\rightarrow\infty} \int_{B_R(0)} \frac{1}{|x-y|^{n-2s}}
  \operatorname{d}\!y\\
   &=& C(n,s)M_R \int_{B_R(0)}\lim_{|x|\rightarrow\infty} \frac{1}{|x-y|^{n-2s}}
  \operatorname{d}\!y=0
\end{eqnarray*}
uniformly for $t\in\mathbb{R}$. Using the same argument as above, we also obtain
\begin{eqnarray*}
 \lim_{|t|\rightarrow\infty} v_R(x,t)
   &\leq& A_{n,s}M_R\lim_{|t|\rightarrow\infty} \int_{-R^2}^{R^2} \int_{B_R(0)} \frac{e^{-\frac{|x-y|^2}{4(t-\tau)}}}{(t-\tau)^{\frac{n}{2}+1-s}}\chi_{\{t>\tau\}}\operatorname{d}\!y\operatorname{d}\!\tau
\\
    &\leq& C(n,s)M_R\lim_{|t|\rightarrow\infty} \int_{-R^2}^{R^2} \frac{1}{|t-\tau|^{1-s}}
  \operatorname{d}\!\tau\\
   &=& C(n,s)M_R \int_{-R^2}^{R^2}\lim_{|t|\rightarrow\infty} \frac{1}{|t-\tau|^{1-s}}
  \operatorname{d}\!\tau=0
\end{eqnarray*}
uniformly for $x\in\mathbb{R}^n$. By virtue of the nonnegativity of $v_R$, thus we conclude that both \eqref{v_Rd1} and \eqref{v_Rd2} hold. Thereby if \eqref{mp} is not valid, then \eqref{w_Rd1} and \eqref{w_Rd2} imply that there must exist $(x^0,t_0)\in \mathbb{R}^n\times\mathbb{R}$ such that
\begin{equation*}
  w_R(x^0,t_0)=\inf_{\mathbb{R}^n\times\mathbb{R}} w_R(x,t)<0.
\end{equation*}
By a direct calculation, we have
\begin{equation*}
  (\partial_t-\Delta)^s w_R(x^0,t_0)
=C_{n,s}\int_{-\infty}^{t}\int_{\mathbb{R}^n}
  \frac{w_R(x^0,t_0)-w_R(y,\tau)}{(t_0-\tau)^{\frac{n}{2}+1+s}}e^{-\frac{|x^0-y|^2}{4(t_0-\tau)}}
  \operatorname{d}\!y\operatorname{d}\!\tau<0,
\end{equation*}
which contradicts the equation \eqref{w_R}. Thus, we verify that \eqref{mp} is true.
Substituting the definition of $v_R$ into \eqref{mp} and letting $R\rightarrow\infty$, then it follows from Levi's monotone convergence theorem that
\begin{eqnarray}\label{uv}
  u(x,t) &\geq& v_R(x,t)\nonumber \\
  &=&  \int_{-\infty}^{+\infty} \int_{\mathbb{R}^n} G(x-y,t-\tau)f_R(y,\tau,u(y,\tau))\operatorname{d}\!y\operatorname{d}\!\tau\nonumber\\
   &\rightarrow&  \int_{-\infty}^{+\infty} \int_{\mathbb{R}^n} G(x-y,t-\tau)f(y,\tau,u(y,\tau))\operatorname{d}\!y\operatorname{d}\!\tau=:v(x,t),
\end{eqnarray}
as $R\rightarrow\infty$ for every $(x,t)\in  \mathbb{R}^n\times\mathbb{R} $.
As a consequence of the above relation that shows
the existence of the integral, if
we denote
\begin{equation*}
  w(x,t):=u(x,t)-v(x,t),
\end{equation*}
then it is plain that $w$ is a solution of
\begin{equation*}
\left\{
\begin{array}{ll}
    (\partial_t-\Delta)^s w(x,t)=0 ,~   &(x,t) \in  \mathbb{R}^n\times\mathbb{R}  , \\
  w(x,t)\geq 0 , ~ &(x,t)  \in \mathbb{R}^n\times\mathbb{R}.
\end{array}
\right.
\end{equation*}
Applying the Liouville theorem established in Theorem \ref{Liouville}\,, we immediately obtain
\begin{equation*}
  w(x,t)\equiv C\geq 0.
\end{equation*}
If $C>0$, then
$$u(x,t)=v(x,t)+C\geq C>0.$$
A combination of the condition \eqref{fcondition} and \eqref{uv} yields that
\begin{eqnarray*}
  u(x,t) &\geq & A_{n,s}C\int_{-\infty}^{t} \int_{\mathbb{R}^n} \frac{e^{-\frac{|x-y|^2}{4(t-\tau)}}}{(t-\tau)^{\frac{n}{2}+1-s}}
  \operatorname{d}\!y\operatorname{d}\!\tau \\
   &=& C(n,s) \int_{-\infty}^{t} \frac{1}{(t-\tau)^{1-s}}\operatorname{d}\!\tau =\infty.
\end{eqnarray*}
This contradiction infers that $C=0$, and thus $u$ satisfies the integral equation
\begin{equation*}
  u(x,t)=\int_{-\infty}^{+\infty} \int_{\mathbb{R}^n} G(x-y,t-\tau)f(y,\tau,u(y,\tau))\operatorname{d}\!y\operatorname{d}\!\tau.
\end{equation*}

On the other hand, if $u(x,t)$ is a solution of the integral equation \eqref{ieq}, then
combining Definition \ref{dlfd} with Fubini's theorem, we derive
\begin{eqnarray*}
  &&\int_{-\infty}^{+\infty} \int_{\mathbb{R}^n} (\partial_t-\Delta)^su(x,t)\varphi(x,t)\operatorname{d}\!x\operatorname{d}\!t\\
  &=&  \int_{-\infty}^{+\infty} \int_{\mathbb{R}^n} u(x,t)(\partial_t-\Delta)_{\rm right}^s\varphi(x,t)\operatorname{d}\!x\operatorname{d}\!t\\
   &=&  \int_{-\infty}^{+\infty} \int_{\mathbb{R}^n} \left(\int_{-\infty}^{+\infty} \int_{\mathbb{R}^n} G(x-y,t-\tau)f(y,\tau,u(y,\tau))\operatorname{d}\!y\operatorname{d}\!\tau\right)
   (\partial_t-\Delta)_{\rm right}^s\varphi(x,t)\operatorname{d}\!x\operatorname{d}\!t\\
   &=&\int_{-\infty}^{+\infty} \int_{\mathbb{R}^n} \left(\int_{-\infty}^{+\infty} \int_{\mathbb{R}^n} G(x-y,t-\tau)(\partial_t-\Delta)_{\rm right}^s\varphi(x,t)\operatorname{d}\!x\operatorname{d}\!t\right)
   f(y,\tau,u(y,\tau))\operatorname{d}\!y\operatorname{d}\!\tau\\
   &=&\int_{-\infty}^{+\infty} \int_{\mathbb{R}^n} f(y,\tau,u(y,\tau))\varphi(y,\tau)\operatorname{d}\!y\operatorname{d}\!\tau\\
      &=&\int_{-\infty}^{+\infty} \int_{\mathbb{R}^n} f(x,t,u(x,t))\varphi(x,t)\operatorname{d}\!x\operatorname{d}\!t
\end{eqnarray*}
for any real test function $\varphi(x,t)\in C_0^\infty(\mathbb{R}^n\times\mathbb{R})$.
This unequivocally indicates that $u$ is also a solution to the pseudo-differential equation \eqref{pdeq}.

As hinted at above, we complete the proof of Theorem \ref{equivalence}\,.
\end{proof}

\section*{Acknowledgments}
The work of the first author is partially supported by MPS Simons foundation 847690, and
the work of the second author is partially supported by the National Natural Science Foundation of China (NSFC Grant No. 12101452).

\end{document}